\documentclass[twoside]{amsart}

\usepackage{latexsym}
\usepackage{amsthm,amsfonts,amssymb}
\usepackage[leqno]{amsmath}
\usepackage{xspace}
\usepackage{cite}
\usepackage[all]{xypic}
\usepackage{hyperref}

\newcommand{\numberseries}{\mdseries}   

\newlength{\thmtopspace}                
\newlength{\thmbotspace}                
\newlength{\thmheadspace}               
\newlength{\thmindent}                  

\newtheoremstyle{bfupright head,slanted body}
                {\thmtopspace}{\thmbotspace}
                {\slshape}{\thmindent}{\bfseries}{.}{\thmheadspace}
                {{\numberseries \thmnumber{(#2) }}\thmnote{#3}}

\newtheoremstyle{bfupright head,upright body}
                {\thmtopspace}{\thmbotspace}
                {\upshape}{\thmindent}{\bfseries}{.}{\thmheadspace}
                {{\numberseries \thmnumber{(#2) }}\thmnote{#3}}

\newtheoremstyle{bfit head,upright body}
                {\thmtopspace}{\thmbotspace}
                {\upshape}{\thmindent}{\upshape}{.}{\thmheadspace}
                {{\numberseries\thmnumber{(#2) }}
                {\bfseries\itshape\thmnote{\negthickspace#3}}}

\newtheoremstyle{it head,upright body}
                {\thmtopspace}{\thmbotspace}
                {\upshape}{\thmindent}{\upshape}{.}{\thmheadspace}
                {{\numberseries\thmnumber{(#2) }}
                {\itshape\thmnote{\negthickspace#3}}}

\newtheoremstyle{fixed bf head,slanted body}
                {\thmtopspace}{\thmbotspace}{\slshape}
                {\thmindent}{\bfseries}{.}{\thmheadspace}
                {{\numberseries \thmnumber{(#2) }}\thmname{#1}\thmnote{ (#3)}}

\newtheoremstyle{fixed bf head,upright body}
                {\thmtopspace}{\thmbotspace}{\upshape}
                {\thmindent}{\bfseries}{.}{\thmheadspace}
                {{\numberseries \thmnumber{(#2) }}\thmname{#1}\thmnote{ (#3)}}

\newtheoremstyle{independent paragraph}
                {\thmtopspace}{\thmbotspace}
                {\upshape}{\thmindent}{\upshape}{}{0pt}
                {\thmnote{#3 }}

\newtheoremstyle{subparagraph}
                {\thmbotspace}{\thmbotspace}
                {\upshape}{\thmindent}{\upshape}{}{0pt}
                {\thmnote{#3 }}

\newtheoremstyle{notes}
                {\thmtopspace}{\thmbotspace}
                {\ttfamily}{\thmindent}{\ttfamily\small }{}{0pt}
                {\thmnote{#3 }}

\theoremstyle{bfupright head,slanted body}
\newtheorem{res}{}[section]             \newtheorem*{res*}{}

\theoremstyle{bfit head,upright body}
                 \newtheorem*{com*}{}

\theoremstyle{bfupright head,upright body}
\newtheorem{bfhpg}[res]{}               \newtheorem*{bfhpg*}{}

\theoremstyle{it head,upright body}
               \newtheorem*{ithpg*}{}

\theoremstyle{fixed bf head,slanted body}
\newtheorem{thm}[res]{Theorem}          \newtheorem*{thm*}{Theorem}
\newtheorem{prp}[res]{Proposition}      \newtheorem*{prp*}{Proposition}
\newtheorem{cor}[res]{Corollary}        \newtheorem*{cor*}{Corollary}
\newtheorem{lem}[res]{Lemma}            \newtheorem*{lem*}{Lemma}

\theoremstyle{fixed bf head,upright body}
       \newtheorem*{dfn*}{Definition}
      \newtheorem*{obs*}{Observation}
\newtheorem{rmk}[res]{Remark}           \newtheorem*{rmk*}{Remark}
\newtheorem{exa}[res]{Example}          \newtheorem*{exa*}{Example}

\theoremstyle{independent paragraph}
\newtheorem{ipg}{}

\theoremstyle{subparagraph}
\newtheorem{spg}{}

\newlength{\thmlistleft}        
\newlength{\thmlistright}       
\newlength{\thmlistpartopsep}   
\newlength{\thmlisttopsep}      
\newlength{\thmlistparsep}      
\newlength{\thmlistitemsep}     

\newcounter{eqc} 
  {\end{list}}%

\newcounter{prt}
\newenvironment{prt}{\begin{list}{\upshape (\alph{prt})}%
    {\usecounter{prt}%
      \setlength{\leftmargin}{\thmlistleft}%
      \setlength{\labelwidth}{\thmlistleft}%
      \setlength{\rightmargin}{\thmlistright}%
      \setlength{\partopsep}{\thmlistpartopsep}%
      \setlength{\topsep}{\thmlisttopsep}%
      \setlength{\parsep}{\thmlistparsep}%
      \setlength{\itemsep}{\thmlistitemsep}}}%
  {\end{list}}%

\newcommand{\prtlbl}[1]{{\upshape(#1)}}

\newcounter{rqm}
  {\end{list}}%

  {\end{list}}%

  {\end{list}}%

\newenvironment{prf}[1][Proof]{\begin{proof}[\bf #1]}{\end{proof}}

  \newcommand{\step}[1]{${#1}^\circ$}

  \newcommand{\proofoftag}[2][:]{(#2)#1}

\newcommand{\pgref}[1]{(\ref{#1})}
\newcommand{\pgpartref}[2]{(\ref{#1})\prtlbl{#2}}

\renewcommand{\eqref}[1]{\pgref{eq:#1}}

\newcommand{\corref}[2][Corollary~]{#1\pgref{cor:#2}}

\newcommand{\exaref}[2][Example~]{#1\pgref{exa:#2}}
\newcommand{\lemref}[2][Lemma~]{#1\pgref{lem:#2}}

\newcommand{\prpref}[2][Proposition~]{#1\pgref{prp:#2}}

\newcommand{\rmkref}[2][Remark~]{#1\pgref{rmk:#2}}
\newcommand{\thmref}[2][Theorem~]{#1\pgref{thm:#2}}

\newcommand{\lempartref}[3][Lemma~]{#1\pgpartref{lem:#2}{#3}}

\newcommand{\corcite}[2][?]{\cite[cor.~#1]{#2}}
\newcommand{\lemcite}[2][?]{\cite[lem.~#1]{#2}}
\newcommand{\prpcite}[2][?]{\cite[prop.~#1]{#2}}
\newcommand{\rescite}[2][?]{\cite[#1]{#2}}

\newcommand{\thmcite}[2][?]{\cite[thm.~#1]{#2}}

\newcommand{\catz}{0}

\newcommand{\Cat}[2]{{\sf{#2}}(#1)}
\newcommand{\Catsup}[3]{{\sf{#2}}^{\text{\upshape #3}}(#1)}
\newcommand{\Catsub}[3]{{\sf{#2}}_{#3}(#1)}

\newcommand{\Catz}[2]{\Catsub{#1}{#2}{\catz}}

\newcommand{\Fz}[1][R]{\Catz{#1}{F}}

\newcommand{\A}[1][R]{\Cat{#1}{A}}

\renewcommand{\H}[2][\no]{\operatorname{H}_{#1}(#2)}
\newcommand{\Susp}[2]{{\scriptstyle\Sigma}^{#1}{#2}}

\newcommand{\Hom}[3][R]{\operatorname{Hom}_{#1}(#2,#3)}
\newcommand{\Ext}[4][R]{\operatorname{Ext}_{#1}^{#2}(#3,#4)}
\newcommand{\tp}[3][R]{#2\otimes_{#1}#3}

\newcommand{\Tor}[4][R]{\operatorname{Tor}^{#1}_{#2}(#3,#4)}

\newcommand{\D}[1][R]{\Cat{#1}{D}}

\newcommand{\Df}[1][R]{\Catsup{#1}{D}{f}}

\newcommand{\DHom}[3][R]{\operatorname{\mathbf{R}Hom}_{#1}(#2,#3)}
\newcommand{\Dtp}[3][R]{#2\otimes_{#1}^{\mathbf{L}}#3}

\newcommand{\no}{\mspace{-1mu}} 

\renewcommand{\a}{\alpha}
\renewcommand{\b}{\beta}
\newcommand{\f}{\varphi}
\newcommand{\g}{\gamma}
\renewcommand{\l}{\ell}

\newcommand{\eq}{\simeq}
\newcommand{\is}{\cong}

\newcommand{\ZZ}{\mathbb{Z}}

\newcommand{\m}{\mathfrak{m}}
\newcommand{\n}{\mathfrak{n}}
\newcommand{\p}{\mathfrak{p}}
\newcommand{\q}{\mathfrak{q}}

\newcommand{\Rmk}{(R,\m,k)}
\newcommand{\Snl}{(S,\n,l)}

\newcommand{\Rhat}{\widehat{R}}
\newcommand{\Shat}{\widehat{S}}

\newcommand{\Lra}{\iff}

\newcommand{\into}{\hookrightarrow}

\newcommand{\xra}{\xrightarrow}

\newcommand{\mapdef}[4][\rightarrow]{\mbox{\ensuremath{#2\!: #3 #1 #4}}}

\newcommand{\supremum}[2]{\sup{\{\,#1\:|\:#2\}}}

\newcommand{\E}[2][R]{\operatorname{E}_{#1}(#2)}

\newcommand{\Ker}[1]{\mbox{\ensuremath{\operatorname{Ker}#1}}}
\newcommand{\Coker}[1]{\mbox{\ensuremath{\operatorname{Coker}#1}}}

\newcommand{\SpecR}{\operatorname{Spec}R}

\newcommand{\dptR}{\operatorname{depth}R}

\newcommand{\dimR}{\operatorname{dim}R}

\newcommand{\Spec}[1]{\operatorname{Spec}#1}
\newcommand{\Max}[2][R]{\operatorname{Max}_{#1}#2}

\newcommand{\Ass}[2][R]{\operatorname{Ass}_{#1}#2}
\newcommand{\Supp}[2][R]{\operatorname{Supp}_{#1}#2}
\newcommand{\supp}[2][R]{\operatorname{supp}_{#1}#2}

\newcommand{\fd}[2][R]{\operatorname{fd}_{#1}#2}

\newcommand{\pd}[2][R]{\operatorname{pd}_{#1}#2}

\newcommand{\Rfd}[2][R]{\operatorname{\mathsf{R}fd}_{#1}#2}

\newcommand{\Gdim}[2][R]{\operatorname{G--dim}_{#1}#2}
\newcommand{\Gfd}[2][R]{\operatorname{Gfd}_{#1}#2}

\newcommand{\dpt}[2][R]{\operatorname{depth}_{#1}#2}

\newcommand{\supP}[1]{\sup{(#1)}}

\newcommand{\tpP}[3][R]{(\tp[#1]{#2}{#3})}

\newcommand{\DtpP}[3][R]{(\Dtp[#1]{#2}{#3})}

\newcommand{\one}{\ensuremath{\mathord \dagger}}
\newcommand{\two}{\ensuremath{\mathord \ddagger}}
\newcommand{\three}{\ensuremath{\mathord \ast\ast}}

\newcommand{\five}{\ensuremath{\mathord \ast}}

\newcommand{\seven}{\ensuremath{\mathord \dagger\dagger}}

\renewcommand{\theequation}{\arabic{equation}}
\numberwithin{equation}{res}

\newcommand{\msupp}[2][R]{\operatorname{max}_{#1}#2}

\begin{document}

\title[Gorenstein dimension of modules over homomorphisms]{Gorenstein
  dimension\\ of modules over homomorphisms}

\author{Lars Winther Christensen\ \ and\ \ Srikanth Iyengar}

\thanks{L.W.C.\ was partly supported by a grant from the Danish
  Natural Science Research Council.}

\thanks{S.I.\ was partly supported by NSF grant DMS 0442242.}

\address{Lars Winther Christensen, Department of Mathematics,
  University of Nebraska, Lincoln, NE 68588-0130, U.S.A.}

\email{winther@math.unl.edu}

\address{Srikanth Iyengar, Department of Mathematics, University of
  Nebraska, Lincoln, NE 68588-0130, U.S.A.}

\email{iyengar@math.unl.edu}

\date{\today}



\begin{abstract}
  Given a homomorphism of commutative noetherian rings $R\to S$ and an
  $S$--module $N$, it is proved that the Gorenstein flat dimension of
  $N$ over $R$, when finite, may be computed locally over $S$.  When,
  in addition, the homomorphism is local and $N$ is finitely generated
  over $S$, the Gorenstein flat dimension equals
  $\supremum{m\in\ZZ}{\Tor{m}{E}{N}\ne 0}$, where $E$ is the injective
  hull of the residue field of $R$. This result is analogous to a
  theorem of Andr\'{e} on flat dimension.
\end{abstract}

\maketitle


\section*{Introduction}

Let $R$ be a commutative noetherian ring and let $N$ be an
$R$--module.  We say that $N$ is \emph{finite over a homomorphism} if
there exists a homomorphism of rings $R\to S$ such that $S$ is
noetherian, $N$ is a finite (that is, finitely generated) $S$--module,
and the $S$--action is compatible with the action of $R$.

In the case where $R \to S$ is a local homomorphism, this class of
modules has been studied by Apassov \cite{DAp99b}, who called them
almost finite modules, and by Avramov, Foxby, Miller, Sather-Wagstaff
and others, cf.~\cite{LLAHBF91,SInSSW04,AIM-}.  The work of these and
other authors show that modules finite over (local) homomorphisms have
homological properties extending those of finite modules (over local
rings).

An important property of many invariants of $R$--modules is that they
can be computed locally over $R$. A basic question is whether the same
property holds for modules over a homomorphism; that is, whether an
invariant of the $R$--module $N$ can be computed locally over $S$.  It
is easy to see that this is the case for flat dimension; this paper
focuses on the Gorenstein flat dimension.  Introduced by Enochs, Jenda
and Torrecillas~\cite{EJT-93}, this invariant is one generalization to
non-finite modules of the notion of G--dimension, due to Auslander and
Bridger \cite{MAs67,MAsMBr69}. In \thmref{loc} we prove that if
$\Gfd{N}$, the Gorenstein flat dimension of $N$, is finite, then
\begin{equation*}
  \Gfd{N} =
  \supremum{\Gfd[R_\p]{N_\q}}{\q\in\Spec{S}\;\text{and}\;\p=R\cap\q}.
\end{equation*}
This extends a well-known result \cite{CFF-02,HHl04a} for the absolute case
$R \xra{=} S$.

The result above focuses attention on modules over local
homomorphisms.  In this situation, a theorem of Andr\'{e} \cite{hac}
says that if $N$ is finite, then the flat dimension over $R$ equals
$\supremum{m\in\ZZ}{\Tor{m}{k}{N}\ne 0}$, where $k$ is the residue
field of $R$.  \thmref{main} gives an analogous result in the context
of Gorenstein flat dimension: If $N$ is finite over a local
homomorphism, and $\Gfd{N}$ is finite, then
\begin{equation*}
  \Gfd{N} = \supremum{m\in\ZZ}{\Tor{m}{E}{N}\ne 0},
\end{equation*}
where $E$ is the injective hull of the residue field $k$. The absolute
case appears in \cite{LWC}.

A crucial difference between this result and Andr\'{e}'s is that it
must be assumed \textsl{a priori} that $\Gfd{N}$ is finite: Vanishing
of $\Tor{\gg 0}{E}{N}$ does not detect finite Gorenstein flat
dimension, see \exaref{js}. This example also suggests that
Andr\'{e}'s proof, which relies on the fact that finite flat dimension
of $N$ is detected by vanishing of $\Tor{\gg 0}{-}{N}$, is not likely
to carry over to our context. And, indeed, our arguments have a
different flavor.

As a corollary we obtain the following result about completions: If
$N$ is finite over a local homomorphism, and $\Gfd{N}$ is finite, then
\begin{equation*}
  \Gfd{N} = \Gfd[\Rhat]{\tpP[S]{\Shat}{N}}.
\end{equation*}
The corresponding result for flat dimension is elementary; for
Gorenstein flat dimension we are not aware of any other proof.


\section{Basic notions}

Throughout the paper $R$ and $S$ denote rings; unless stated
otherwise, they are assumed to be commutative and noetherian. Given a
homomorphism $\mapdef{\f}{R}{S}$, any $S$--module becomes an
$R$--module with the action determined by $\f$. We say that $\f$ is
\emph{local}, if $R$ and $S$ are local rings with maximal ideals $\m$
and $\n$, and $\f(\m) \subseteq \n$.

We work with complexes, which we grade homologically:
\begin{equation*}
  M = \dots \to M_{\l+1} \to M_\l \to M_{\l-1} \to \cdots
\end{equation*}
The homological size of a complex is captured by the numbers $\sup{M}$
and $\inf{M}$, defined as the supremum and infimum of the set
$\{\l\in\ZZ \mid \H[\l]{M}\ne 0\}$. We say that $M$ is
\emph{homologically finite} if the $R$--module $\H{M}$ is
\emph{finite}, that is, finitely generated.

We use the notation $\D$ for the derived category of $R$, and $\Df$
for its subcategory of homologically finite complexes.  We use the
symbol $\eq$ to denote isomorphisms in derived categories.

Let $L$ and $M$ be $R$--complexes, that is to say, complexes of
$R$--modules.  The derived tensor product and Hom functors are denoted
$\Dtp{L}{M}$ and $\DHom{L}{M}$. We write $\pd{M}$ for the projective
dimension, and $\fd{M}$ for the flat dimension, of $M$ over $R$,
cf.~\cite{LLAHBF91}.

When $\Rmk$ is local, the \emph{depth} of an $R$--complex $M$ is
defined by
\begin{equation}
  \label{eq:dpt}
  \dpt{M} = -\sup{\DHom{k}{M}}.
\end{equation}
Thus, $\dpt{M}=\infty$ if and only if $\H{\DHom{k}{M}}=0$.

\begin{bfhpg}[Supports]
  \label{supp}
  The \emph{support} of an $R$--complex $M$ is a subset of $\SpecR$:
  \[
  \Supp{M}=\{\p\in\SpecR \mid \H{M_\p}\ne 0\},
  \]
  and $\Max{M}$ is the subset of maximal ideals in $\Supp{M}$.
  
  Foxby \cite{HBF79} has introduced the \emph{small support} of $M$ as
  \[
  \supp{M} =\{\p\in \SpecR\mid \H{\Dtp{k(\p)}{M}}\ne 0\},
  \] 
  where $k(\p)$ denotes the residue field $R_\p/\p R_\p$ of $R$ at
  $\p$.  For convenience we set
  \[
  \msupp{M} =\{\p\in\supp{M}\mid \text{$\p$ is maximal in
    $\supp{M}$}\}.
  \]
  
  Note that $\supp{M} \subseteq \Supp{M}$ and equality holds when $M$
  is homologically finite. Elements in $\Max{M}$ are maximal ideals,
  while those in $\msupp{M}$ need not be; see property \prtlbl{c}
  below.

  \begin{spg}
    We recall some properties of these subsets. Let $L$ and $M$ be
    $R$--complexes. For $\p$ in $\SpecR$ we write $\E{R/\p}$ for the
    injective hull of the $R$--module $R/\p$.
    \begin{prt}
    \item $\supp{M} = \emptyset$ if and only if $\H{M}=0$.
    \item $\supp{(\Dtp{L}{M})}=\supp{L}\cap\supp{M}$, for any
      $R$--complex $L$.
    \item $\supp{\E{R/\p}} = \{\p\}$, for any $\p$ in $\SpecR$.
    \item A prime ideal $\p$ is in $\supp{M}$ if and only if
      $\H{\Dtp{\E{R/\p}}{M}}\ne 0$.
    \item When $\sup{M}=s$ is finite, the associated primes of the top
      homology module belong to the small support:
      $\Ass{\H[s]{M}}\subseteq \supp{M}$.
    \item When $\Rmk$ is local, $\m$ is in $\supp{M}$ if and only if
      $\dpt{M}$ is finite.
    \end{prt}
    
    Indeed, parts \prtlbl{a}, \prtlbl{b}, \prtlbl{c}, and \prtlbl{f}
    are proved in \cite[sec.~2]{HBF79}; part \prtlbl{d} follows
    immediately from \prtlbl{a}, \prtlbl{b} and \prtlbl{c}; part
    \prtlbl{e} is \prpcite[2.6 and (2.4.1)]{LWC01b}.
  \end{spg}
\end{bfhpg}

\begin{ipg}
  Next we recall the notion of G--dimension; see
  \cite{MAsMBr69,LLAAMr02,LWC} for details.
\end{ipg}

\begin{bfhpg}[G--dimension]
  \label{dfn:gdim}
  A finite $R$--module $G$ is said to be \emph{totally reflexive} if
  there exists an exact complex $L$ of finite free $R$--modules such
  that $G \is \Coker{(L_1 \to L_0)}$ and $\H{\Hom{L}{R}}=0$. Any
  finite free module is totally reflexive, so each homologically
  finite $R$--complex $N$ with $\H[\l]{N}=0$ for $\l\ll 0$ admits a
  resolution by totally reflexive modules. The \emph{G--dimension} is
  the number
  \begin{gather*}
    \Gdim{N} = \inf\left\{d \in\ZZ \left|
        \begin{gathered}
          \text{$N$ is isomorphic in $\D$ to a complex of totally}\\
          \text{reflexive modules: } 0 \to G_d \to G_{d-1} \to \dots
          \to G_i \to 0
        \end{gathered}
      \right\}\right.
  \end{gather*}
\end{bfhpg}

\begin{ipg}
  Enochs, Jenda and Torrecillas~\cite{EJT-93,EEnOJn95b} have studied
  extensions of G--dimension to complexes whose homology may not be
  finite. One such extension is the Gorenstein flat dimension; see
  \cite{EJT-93,LWC}.
\end{ipg}

\begin{bfhpg}[Gorenstein flat dimension]
  An $R$--module $A$ is \emph{Gorenstein flat} if there exists an
  exact complex $F$ of flat modules such that $A \is \Coker{(F_1 \to
    F_0)}$ and $\H{\tp{J}{F}}=0$ for any injective $R$--module $J$.
  Any free module is Gorenstein flat, so each complex $M$ with
  $\H[\l]{M}=0$ for $\l\ll 0$ admits a resolution by Gorenstein flat
  modules. The \emph{Gorenstein flat dimension} is the number
  \begin{gather*}
    \Gfd{M} = \inf\left\{d \in\ZZ \left|
        \begin{gathered}
          \text{$M$ is isomorphic in $\D$ to a complex of Gorenstein}\\
          \text{flat modules: }0 \to A_d \to A_{d-1}\to \dots \to A_i
          \to 0
        \end{gathered}
      \right\}\right.
  \end{gather*}
  When $M$ is homologically finite $\Gfd{M} = \Gdim{M}$;
  see~\thmcite[(5.1.11)]{LWC}.
\end{bfhpg}

\begin{rmk}
  \label{rmk:holm}
  By \thmcite[(3.5) and cor.~(3.6)]{CFH-}, if $M$ is an $R$--complex
  of finite Gorenstein flat dimension, then:
  \begin{align*}
    \Gfd{M} &= \supremum{\supP{\Dtp{J}{M}}}{J\ \text{is injective}}\\
    &= \supremum{\supP{\Dtp{\E{R/\p}}{M}}}{\p\in\SpecR}.
  \end{align*}
\end{rmk}


\section{Localization}
\label{sec:loc}

The gist of this section is that for complexes over homomorphisms the
Gorenstein flat dimension, when it is finite, may be computed locally.
We should like to note that the analogue for flat dimensions is
elementary to verify, for the finiteness of that invariant is detected
by vanishing of Tor functors. The absolute case, $R \xra{=} S$, is
easily deduced from \thmcite[8.8]{SInSSW04} and \thmcite[(2.4)]{CFF-02}.
\begin{thm}
  \label{thm:loc}
  Let $\mapdef{\f}{R}{S}$ be a homomorphism of rings and let $X$ be an
  $S$--complex. If $\,\Gfd{X}$ is finite, then
  \begin{align*}
    \Gfd{X} =\begin{cases} \supremum{\Gfd[R_\p]{X_\q}}{\q\in
        \Spec{S}\;\text{\rm
          and}\;\,\p=\q\cap R}\\
      \supremum{\Gfd[R_\p]{X_\q}}{\q\in \Max[S]{X}\;\text{\rm
          and}\;\,\p=\q\cap R}\\
      \supremum{\Gfd[R_\p]{X_\q}}{\q\in \msupp[S]{X}\;\text{\rm
          and}\;\,\p=\q\cap R}
    \end{cases}
  \end{align*}
\end{thm}

\begin{ipg}
  The proof is given towards the end of this section. In preparation
  we recall a result about colimits of Gorenstein flat modules:
\end{ipg}

\begin{rmk}
  \label{rmk:limGflat}
  If $(M_i)_{i\in I}$ is a filtered system of Gorenstein flat modules
  over a coherent ring, then the colimit $\varinjlim{M_i}$ is
  Gorenstein flat. This follows from work of Enochs et.~al.\ 
  \cite{EEnJLR02,EJO-01} and Holm \cite{HHl04a}: By \thmcite[2.4 (and
  remarks before sec.~2)]{EEnJLR02} a filtered colimit $M =
  \varinjlim{M_i}$ of Gorenstein flat modules has a co-proper right
  resolution by flat modules. Because colimits commute with tensor
  products, \rmkref[]{holm} provides an equality
  \begin{equation*}
    \supremum{\supP{\Dtp{J}{M}}}{J\ \text{is injective}} =0.
  \end{equation*}
  Therefore, by \thmcite[3.6]{HHl04a}, the colimit $M$ is Gorenstein
  flat.
\end{rmk}

\begin{ipg}
  For the next result note that any $R$--module $M$ has a natural
  structure of a module over its endomorphism ring $\Hom{M}{M}$.
\end{ipg}

\begin{lem}
  \label{lem:endo}
  Let $R$ be a coherent ring and $M$ a Gorenstein flat $R$--module.
  Let $Z$ be a multiplicatively closed set in the center of the ring
  $\Hom{M}{M}$.  Then the $R$--module $Z^{-1}M$ is Gorenstein flat.
\end{lem}

\begin{prf}
  Let $\mathcal{V}$ denote the set of finitely generated (as
  semigroups) multiplicatively closed subsets of $Z$. The modules
  $V^{-1}M$, for $V\in\mathcal{V}$, with natural maps
  \begin{equation*}
    \mapdef{\rho^{UV}}{U^{-1}M}{V^{-1}M}\quad\text{for }U\subseteq V
  \end{equation*}
  form a filtered system. It is straightforward to verify that the
  colimit $\varinjlim{V^{-1}M}$ is isomorphic to $Z^{-1}M$ as
  $\Hom{M}{M}$--module and, therefore, as an $R$--module.
  
  By \rmkref{limGflat}, a filtered colimit of Gorenstein flat modules
  is Gorenstein flat, so it remains to see that the modules $V^{-1}M$
  are Gorenstein flat. For any $V\in\mathcal{V}$ the module $V^{-1}M$
  can be constructed by successively inverting the finitely many
  generators of $V$. Thus, it suffices to prove that $M_z$ is
  Gorenstein flat for any $z\in Z$. Again, $M_z$ is the colimit of the
  linear system $(M \xra{z} M \xra{z} M \xra{z} \cdots)$ and hence
  Gorenstein flat by \rmkref[]{limGflat}.
\end{prf}

\begin{ipg}
  We should like to stress that in the next result the ring $S$ need
  not be noetherian.
\end{ipg}

\begin{prp}
  \label{prp:loc}
  Let $R$ be a noetherian ring. Let $\mapdef{\f}{R}{S}$ be a
  homomorphism of rings and $X$ an $S$--complex. For each
  $\q\in\Spec{S}$ and $\p=\q\cap R$, one has
  \begin{equation*}
    \Gfd[R_\p]{X_\q} = \Gfd{X_\q} \le \Gfd{X}.
  \end{equation*}
\end{prp}

\begin{prf}
  The equality in the statement is evident: a Gorenstein flat
  $R_\p$--module is Gorenstein flat over $R$ and any Gorenstein flat
  $R$--module localizes to give a Gorenstein flat $R_\p$--module.
  
  In verifying the inequality one may assume that $\Gfd{X}$ is finite.
  Pick a surjective homomorphism $\tilde S\to S$ where $\tilde S$ is
  an $R$--algebra, free as an $R$--module, and let $\tilde\q$ be the
  preimage of $\q$ in $\tilde S$. Evidently, $X_{\tilde\q}\simeq X_\q$
  as $\tilde S$--complexes, and hence also as $R$--complexes, so
  replacing $S$ with $\tilde S$, we assume henceforth that the
  $R$--module $S$ is free.
  
  Let $U$ be a free resolution of $X$ over $S$ and set $\Omega =
  \Ker(\partial_{d-1}^U)$ for $d=\Gfd{X}$. Since $S$ is a free
  $R$--module, $U$ is also an $R$--free resolution of $X$, and since
  $\Gfd{X}$ is finite, $\Omega$ viewed as an $R$--module is Gorenstein
  flat.  Note that one has isomorphisms
  \[
  U_{\q}\simeq X_\q\quad\text{and}\quad
  \Ker(\partial_{d-1}^{U_{\q}})\cong \Omega_{\q}.
  \]
  The complex $U_\q$ consists of flat $R$--modules, so to settle the
  claim it suffices to prove that the $R$--module $\Omega_\q$ is
  Gorenstein flat.  Therefore, it suffices to verify the result in the
  case where the $S$--module $X$ is Gorenstein flat over $R$.
  
  Homothety provides a homomorphism of rings $S\to \Hom{X}{X}$. Let
  $Z$ be the image of $S\setminus\q$ under this map; it is a
  multiplicatively closed subset in the center of $\Hom{X}{X}$, and
  $Z^{-1}X\cong X_\q$ as $S$--modules. It now remains to invoke
  \lemref{endo}.
\end{prf}

\vspace*{0pt}
\begin{prf}[Proof of \thmref{loc}]
  \label{prf:loc}
  \prpref{loc} implies the first inequality below
  \begin{align*}
    \Gfd{X} &\ge \supremum{\Gfd[R_\p]{X_\q}}{\q\in \Spec{S}\;\text{\rm
        and}\;\,\p=\q\cap R}\\
    &\ge \supremum{\Gfd[R_\p]{X_\q}}{\q\in \Max[S]{X}\;\text{\rm
        and}\;\,\p=\q\cap R}\\
    &\ge \supremum{\Gfd[R_\p]{X_\q}}{\q\in \msupp[S]{X}\;\text{\rm
        and}\;\,\p=\q\cap R}.
  \end{align*}
  The second inequality holds because of the inclusion $\Max[S]{X}
  \subseteq \Supp[S]{X}$, and the third follows also by \prpref{loc}
  as any ideal in $\msupp[S]{X}$ is contained in an ideal from
  $\Max[S]{X}$.  This leaves us one inequality to verify:
  \begin{equation*}
    \Gfd{X} \le \supremum{\Gfd[R_\p]{X_\q}}{\q\in \msupp[S]{X}\;\text{\rm
        and}\;\,\p=\q\cap R}.
  \end{equation*}
  Set $d=\Gfd{X}$ and pick a $\tilde{\p}$ in $\SpecR$ for which
  $\supP{\Dtp{\E{R/\tilde{\p}}}{X}}=d$. Pick a prime ideal $\q'$
  associated to the $S$--module $\H[d]{\Dtp{\E{R/\tilde{\p}}}{X}}$ and
  set $\p'=\q'\cap R$. In the (in)equalities below:
  \begin{equation*}
    \begin{split}
      d &=\supP{\Dtp{\E{R/\tilde{\p}}}{X}}\\
      &= \supP{\Dtp{\E{R/\tilde{\p}}}{X_{\q'}}}\\
      &= \supP{\Dtp[R_{\p'}]{\E[R_{\p'}]{R_{\p'}/\tilde{\p}R_{\p'}}}{X_{\q'}}}\\
      &\le \Gfd[R_{\p'}]{X_{\q'}}
    \end{split}    \tag{\one}
  \end{equation*}
  the second one holds by choice of $\q'$, while the third holds
  because $\E{R/\tilde{\p}}$ is an $R_{\p'}$--module, as $\tilde{\p}
  \subseteq \p'$. By \pgpartref{supp}{e} the ideal $\q'$ is in the
  small support of the $S$--complex $\Dtp{\E{R/\tilde{\p}}}{X}$.  The
  first equality below is due to the associativity of the tensor
  product
  \begin{align*}
    \supp[S]{\DtpP{\E{R/\tilde{\p}}}{X}} &=
    \supp[S]{\DtpP[S]{\DtpP{\E{R/\tilde{\p}}}{S}}{X}}\\
    &= \supp[S]{\DtpP{\E{R/\tilde{\p}}}{S}} \cap \supp[S]{X}
  \end{align*}
  while the second one is by \pgpartref{supp}{b}. These show that
  $\q'$ is in $\supp[S]{X}$. Finally, choose $\q\in\msupp[S]{X}$
  containing $\q'$ and set $\p=\q\cap R$. It follows by (\one) and
  \prpref{loc} that $d \le \Gfd[R_{\p'}]{X_{\q'}} \le
  \Gfd[R_{\p}]{X_{\q}}$.
\end{prf}

\section{Approximations}
\label{sec:approx}

In this section we establish an approximation theorem for complexes of
finite G--dimension; this is an important ingredient in the proof of
\thmref{main}. It is a common generalization to complexes of
\thmcite[2.10]{HHl04a} and \lemcite[(2.17)]{CFH-}, which deal with
modules. Similar extensions have been obtained by Holm et.~al.
\cite{FHS-,HHlPJr2}; see \rmkref[]{inf} and the remarks following the
statement of the theorem for further relations to earlier work.

\begin{thm}
  \label{thm:fgapprox}
  Let $S$ be a ring and $N$ a homologically finite $S$--complex with
  finite G--dimension. For each integer $n \le \Gdim[S]{N}$ there
  exists an exact triangle
  \begin{equation*}
    N \to P \to H \to \Susp{}{N}
  \end{equation*}
  in $\Df[S]$ with the following properties:
  \begin{prt}
  \item $\pd[S]{P} = \Gdim[S]{N}$ and $\Gdim[S]{H} \le n$.
  \item There are inequalities: $\inf{P} \ge n \ge \sup{H}$, and
    \begin{equation*}
      \max{\{n,\sup{N}\}} \ge \sup{P}\quad \text{and}\quad \inf{H} \ge
      \min{\{n,\inf{N}+1\}}.
    \end{equation*}
    Moreover, the following induced sequence of $S$--modules is exact:
    \begin{equation*}
      0 \to \H[n]{N} \to \H[n]{P} \to \H[n]{H} \to \H[n-1]{N} \to 0.
    \end{equation*}
  \end{prt}
\end{thm}

We precede the proof with a couple of remarks and a lemma.

\begin{rmk}
  As above, let $N$ be a homologically finite $S$--complex of finite
  G--dimension. By rotating the exact triangle in \thmref[]{fgapprox},
  we see that for each integer $n \le \Gdim[S]{N}$ there exists an
  exact triangle
  \begin{equation*}
    P' \to H' \to N \to \Susp{}{P'}
  \end{equation*}
  in $\Df[S]$ where $\pd[S]{P'} = \Gdim[S]{N}-1$ and $\Gdim[S]{H'} \le
  n-1$.
\end{rmk}

\begin{rmk}
  Let $N$ be a finite $S$--module with finite G--dimension. Applying
  \thmref{fgapprox} with $n=0$ we get from part \prtlbl{b} an exact
  sequence of finite modules
  \begin{equation*}
    0 \to N \to \H[0]{P} \to \H[0]{H} \to 0.
  \end{equation*}
  Moreover, $\H[\l]{H} = 0 = \H[\l]{P}$ for $\l \ne 0$, so from part
  \prtlbl{a} it follows that $\H[0]{H}$ is totally reflexive and
  $\pd{\H[0]{P}} = \Gdim[S]{N}$. Thus we recover
  \lemcite[(2.17)]{CFH-}.
  
  Analogously, if $\Gdim[S]{N}\ge 1$, applying \thmref{fgapprox} with
  $n=1$ yields an exact sequence of finite modules:
  \begin{equation*}
    0 \to \H[1]{P} \to \H[1]{H} \to N \to 0,
  \end{equation*}
  where $\H[1]{H}$ is totally reflexive and $\pd{\H[1]{P}} =
  \Gdim[S]{N}-1$. In this way we also recover \thmcite[2.10]{HHl04a}.
\end{rmk}

\begin{lem}
  \label{lem:PO-sequence}
  Let $X$ be an $S$--complex.  For any injective homomorphism
  $\mapdef{\iota}{X_n}{Y_n}$ of $\,S$--modules there is a commutative
  diagram
  \begin{equation*}
    \xymatrix{X = \cdots \ar[r] & X_{n+1} \ar@{=}[d] \ar[r]^{\a} & X_n
      \ar@{^{(}->}[d]_{\iota} \ar[r]^{\b} & X_{n-1} \ar@{^{(}->}[d]_{\iota'}
      \ar[r]^{\g} & X_{n-2} \ar@{=}[d] \ar[r] & \cdots\\
      Y = \cdots \ar[r] & X_{n+1} \ar[r]^{\a'} & Y_n \ar[r]^{\b'} &
      Y_{n-1} \ar[r]^{\g'} & X_{n-2} \ar[r] & \cdots }
  \end{equation*}
  such that $Y$ is a complex, $\Coker{\iota'} \is \Coker{\iota}$, and
  the induced map $\H{X} \to \H{Y}$ is an isomorphism. When $X_{n-1}$
  and $Y_n$ are finite, $Y_{n-1}$ can be chosen finite.
\end{lem}

\begin{prf}
  Set $\a' = \iota\a$ and let $\mapdef{\b'}{Y_n}{Y_{n-1}}$ be the
  pushout of $\b$ along $\iota$; thus
  \begin{equation*}
    Y_{n-1}= \frac{Y_n\oplus X_{n-1}}{\{(\iota(x),\b(x))\mid x\in X_n\}}.
  \end{equation*}
  Let $\mapdef{\iota'}{X_{n-1}}{Y_{n-1}}$ be the induced map, which
  sends $x$ to $(0,x)$; it is injective because $\iota$ is. Define
  $\mapdef{\g'}{Y_{n-1}}{X_{n-2}}$ by $(y,x) \mapsto \g(x)$. By
  construction the diagram is commutative. It is elementary to check
  that $Y$ is a complex, and the induced map $\Coker{\iota} \to
  \Coker{\iota'}$ an isomorphism. Thus, the cokernel of the inclusion
  of complexes $X \into Y$ is exact, and hence the induced map $\H{X}
  \to \H{Y}$ is bijective. By construction, $Y_{n-1}$ is finite when
  $X_{n-1}$ and $Y_n$ are so.
\end{prf}

\vspace*{0pt}
\begin{prf}[Proof of \thmref{fgapprox}]
  The hypothesis is that $N$ is a homologically finite $S$--complex
  with finite G--dimension; set $d=\Gdim{N}$ and $i=\inf{N}$. Let
  \begin{equation*}
    \dots \to P_{\l} \to P_{\l-1} \to \dots \to P_i \to 0
  \end{equation*}
  be a projective resolution of $N$ by finite modules. For integers
  $n\le d+1$ we construct, by descending induction on $n$, complexes
  $C(n)$ isomorphic to $N$ in $\D[S]$ and of the form
  \begin{equation*}
    C(n)= 0 \to Q_d \to \dots \to Q_{n} \to G_{n-1} \to
    P_{n-2} \to \dots \to P_i \to 0,
  \end{equation*}
  where the modules $Q_\l$ are also finite projective and $G_{n-1}$ is
  totally reflexive.  For the first step, set $G_d = \Ker{(P_{d} \to
    P_{d-1})}$; this module is totally reflexive, the complex $$C(d+1)
  = 0 \to G_d \to P_{d-1} \to \dots \to P_i \to 0$$
  is isomorphic to
  $N$ in $\D[S]$ and has the desired form.  Next we construct $C(n)$
  from $C(n+1)$.  The totally reflexive module $G_n$ in $C(n+1)$
  embeds into a finite free module $\mapdef{\iota}{G_n}{Q_{n}}$ such
  that $\Coker{\iota}$ is totally reflexive.  By \lemref{PO-sequence}
  we have a commutative diagram
  \begin{equation*} 
    \xymatrix{C(n+1) = \cdots \ar[r] & Q_{n+1} \ar@{=}[d]
      \ar[r] & G_n \ar@{^{(}->}[d]_{\iota} \ar[r] & P_{n-1}
      \ar@{^{(}->}[d]_{\iota'} \ar[r] & P_{n-2} \ar@{=}[d] \ar[r] &
      \cdots\\
      C(n)  =  \ar[r] \cdots & Q_{n+1} \ar[r] & Q_{n} \ar[r] & G_{n-1} \ar[r]
      & P_{n-2} \ar[r] & \cdots }
  \end{equation*}
  The module $\Coker{\iota'}$ is isomorphic to $\Coker{\iota}$ and
  hence totally reflexive; therefore $G_{n-1}$ is totally reflexive.
  The complex $C(n)$ has the desired form, by construction, and is
  isomorphic to $C(n+1) \eq N$, again by \lemref[]{PO-sequence}.
  
  Now, fix an integer $n\le d$ and replace $N$ by $C(n)$. Let $P$ be
  the truncation $N_{\geqslant n}$ of $N$ and $H =
  \Susp{}{(N_{\leqslant n-1})}$; the canonical surjection $N\to P$
  yields an exact triangle
  \begin{equation*}
    \tag{\ensuremath{\Delta}}
    N \to P \to H \to \Susp{}{N}.
  \end{equation*}
  We now verify that this triangle has the desired properties:
  
  \proofoftag{a} It is evident from the construction that $\Gdim[S]{H}
  \le n$ and $\pd[S]{P} \le d$. To see that $\pd[S]{P} = d$, apply
  $\DHom[S]{-}{S}$ to ($\Delta$) and take homology to get the exact
  sequence
  \begin{equation*}
    \Ext[S]{d}{P}{S} \to \Ext[S]{d}{N}{S} \to \Ext[S]{d+1}{H}{S}.
  \end{equation*}
  Recall that $\Gdim[S]{X} = \supremum{m\in\ZZ}{\Ext[S]{m}{X}{S} \ne
    0}$ for any homologically finite $S$--complex of finite
  G--dimension, cf.~\corcite[(2.3.8)]{LWC}. Therefore, in the exact
  sequence above, the module on the right is zero as $d \ge n \ge
  \Gdim[S]{H}$, while the middle one is non-zero as $\Gdim[S]{N}=d$.
  Thus, $\Ext[S]{d}{P}{S} \ne 0$.
  
  \proofoftag{b} By construction $\inf{P} \ge n \ge \sup{H}$, so the
  homology exact sequence
  \begin{equation*}
    \dots \to \H[\l]{N} \to \H[\l]{P} \to \H[\l]{H} \to \H[\l-1]{N}
    \to \cdots 
  \end{equation*}
  associated with ($\Delta$) gives the desired exact sequence and
  isomorphisms $\H[\l]{P} \is \H[\l]{N}$ for $\l \ge n+1$ and
  $\H[\l]{H} \is \H[\l-1]{N}$ for $\l \le n-1$. In particular,
  $\max{\{n,\sup{N}\}} \ge \sup{P}$ and $\inf{H} \ge
  \min{\{n,\inf{N}+1\}}$.
\end{prf}

\begin{rmk}
  \label{rmk:inf}
  We note that with G--dimension replaced by Gorenstein projective
  dimension, or by Gorenstein flat dimension, the arguments in the
  preceding proof carry over to the case where the homology modules of
  $N$ are not necessarily finite. In this paper we only need the
  version stated in \thmref{fgapprox}.
\end{rmk}


\section{Local homomorphisms}
\label{sec:smth}

The main result of this section is:

\begin{thm}
  \label{thm:main}
  Let $\Rmk$ be a local ring, and let $N$ be an $R$--complex, finite
  over a local homomorphism. If $\,\Gfd{N}$ is finite, then
  \begin{equation*}
    \Gfd{N}  = \supP{\Dtp{\E{k}}{N}} = \dpt[\no]{R} -\dpt{N}.
  \end{equation*}
\end{thm}

\begin{ipg}
  The second equality was proved in \thmcite[8.7]{SInSSW04}; the
  theorem is motivated by the following considerations:
\end{ipg}

\begin{bfhpg}[Remarks]
  The flat dimension of $N$ can be tested by cyclic modules, $R/\p$,
  and if $\fd{N}$ is finite, then
  \begin{equation*}
    \supP{\Dtp{k}{N}} = \dpt R - \dpt{N}.
  \end{equation*}
  This is the Auslander--Buchsbaum formula for $N$,
  cf.~\cite[p.~153]{HBF79}. Analogously, the Gorenstein flat dimension
  is tested by modules $\E{R/\p}$, cf.~\rmkref[]{holm}, and if
  $\Gfd{N}$ is finite, an analogue of the Auslander--Buchsbaum is
  provided by \thmcite[8.7]{SInSSW04}:
  \begin{equation*}
    \supP{\Dtp{\E{k}}{N}} = \dpt[\no]{R} -\dpt{N}.
  \end{equation*}
  Assume that $N$ is finite over a local homomorphism, then
  \begin{equation*}
    \fd{N} = \supP{\Dtp{k}{N}}
  \end{equation*}
  by \prpcite[5.5]{LLAHBF91}. By the three displayed equations it
  follows that
  \begin{equation}
    \label{eq:fd}
    \fd{N} = \supP{\Dtp{\E{k}}{N}} \;\; \text{when $\fd{N}$ is
      finite}.
  \end{equation}
  However, an elementary argument is also available: Set $f=\fd{N}$;
  associated to the exact sequence $0 \to k \to \E{k} \to C \to 0$ is
  an exact sequence of homology modules
  \begin{equation*}
    0 \to \H[f]{\Dtp{k}{N}} \to \H[f]{\Dtp{\E{k}}{N}} \to \cdots,
  \end{equation*}
  which shows that also $\H[f]{\Dtp{\E{k}}{N}} \ne 0$.
  
  \thmref{main} is an analogue of \eqref{fd} for Gorenstein flat
  dimension. When $N$ is finite over $R$ itself, the first equality in
  \thmref[]{main} recovers \thmcite[(2.4.5)(b)]{LWC}:
  \begin{equation*}
    \Gdim{N} = \supP{\Dtp{\E{k}}{N}}.
  \end{equation*}
  Even in this case one has to assume a priori that the dimension is
  finite:
\end{bfhpg}

\begin{exa}
  \label{exa:js}
  Jorgensen and \c{S}ega \thmcite[1.7]{DAJLMS} construct an artinian
  ring $R$ and a finite $R$--module $L$ with
  \begin{equation*}
    \Gdim{L} = \infty \quad\text{and}\quad \inf{\DHom{L}{R}}=0.
  \end{equation*}
  The last equality translates to $\supP{\Dtp{\E{k}}{L}} = 0$ by
  Matlis duality.
\end{exa}

\begin{ipg}
  It is implicit in \thmref{main} that both $\supP{\Dtp{\E{k}}{N}}$
  and $\dpt{N}$ are finite. This holds in general for complexes finite
  over local homomorphisms:
\end{ipg}

\begin{lem}
  \label{lem:supp}
  Let $\Rmk$ be a local ring and let $N$ be an $R$--complex, finite
  over a local homomorphism. If $\,\H{N} \ne 0$, then
  \begin{equation*}
    \dpt{N} \text{ is finite}\quad \text{and}\quad \H{\Dtp{\E{k}}{N}} \ne 0.
  \end{equation*}
\end{lem}

\begin{prf}
  By assumption there is a local homomorphism
  $\mapdef{\f}{\Rmk}{\Snl}$, such that $N$ is homologically finite
  over $S$. With $i=\inf{N}$ one has
  \begin{equation*}
    \H[i]{\Dtp{k}{N}} \is \tp{k}{\H[i]{N}} \is \H[i]{N}/\m\H[i]{N}.
  \end{equation*}
  Since $\f$ is local, and the $S$--module $\H[i]{N}$ is finite and
  non-zero, Nakayama's lemma implies $\H[i]{N}/\m\H[i]{N}$ is
  non-zero. Thus $\m$ is in $\supp{N}$, in particular, $\dpt{N}$ is
  finite, cf.~\pgpartref{supp}{f}.  Moreover, $\m$ is in
  $\supp{\E{k}}$, by \pgpartref{supp}{c}, and thus also in
  $\supp{\DtpP{\E{k}}{N}}$, whence $\H{\Dtp{\E{k}}{N}} \ne 0$ by
  \pgpartref{supp}{a}.
\end{prf}

\begin{ipg}
  For \thmref{main} it is important that the homology of
  $\Dtp{\E{k}}{N}$ is non-zero. However, that condition alone is not
  sufficient for the first equality, not even for \eqref{fd}; one
  needs the finiteness of $\H{N}$:
\end{ipg}

\begin{exa}
  Let $\Rmk$ be a regular local ring. For a prime ideal $\p\ne \m$ set
  $k(\p) = R_\p/\p R_\p$ and $N= k(\p)\oplus R$. Then
  \begin{gather*}
    \fd{N} = \fd[R_\p]{k(\p)} = \dimR_\p\quad\text{and}\\
    \sup{\DtpP{\E{k}}{N}} = 0= \sup{\DtpP{k}{N}}.
  \end{gather*}
\end{exa}

\begin{ipg}
  For the proof of the theorem we need the following lemmas. The first
  one deals with the restricted flat dimension, introduced by Foxby in
  \cite{CFF-02}. Its relevance for our purpose comes from
  \thmcite[8.8]{SInSSW04}, see also \thmcite[3.19]{HHl04a}.
  
  As usual, for any local ring $\Rmk$ its $\m$-adic completion is
  denoted $\Rhat$.
\end{ipg}
\clearpage

\begin{lem}
  \label{lem:Rfd}
  Let $\mapdef{\f}{R}{S}$ be a homomorphism of rings and $X$ an
  $S$--complex.
  \begin{prt}
  \item If $\,\f$ is flat, then $\Rfd{X} \le \Rfd[S]{X}$
  \item If $\,\f$ is local, then $\Rfd{X} = \Rfd{\tpP[S]{\Shat}{X}} \le
    \Rfd[\Rhat]{\tpP[S]{\Shat}{X}}$
  \end{prt}
\end{lem}

\begin{prf}
  Let $\Fz$ be the class of $R$--modules of finite flat dimension.
  
  \proofoftag{a} For each $T\in\Fz$ the module $\tp{T}{S}$ has finite
  flat dimension over $S$. With this, the desired inequality follows
  from:
  \begin{align*}
    \Rfd{X} &= \supremum{\supP{\Dtp{T}{X}}}{T\in\Fz}\\
    &= \supremum{\supP{\Dtp{T}{\tpP[S]{S}{X}}}}{T\in\Fz}\\
    &= \supremum{\supP{\Dtp[S]{\tpP{T}{S}}{X}}}{T\in\Fz}\\
    &\le \Rfd[S]{X},
  \end{align*}
  where the first equality is the definition.
  
  \proofoftag{b} The inequality follows from \prtlbl{a}, and the
  equality is an easy calculation:
  \begin{align*}
    \Rfd{X} &= \supremum{\supP{\Dtp{T}{X}}}{T\in\Fz}\\
    &= \supremum{\supP{\tp[S]{\DtpP{T}{X}}{\Shat}}}{T\in\Fz}\\
    &= \supremum{\supP{\Dtp{T}{\tpP[S]{X}{\Shat}}}}{T\in\Fz}\\
    &= \Rfd{\tpP[S]{X}{\Shat}}\qedhere
  \end{align*}
\end{prf}

\begin{lem}
  \label{lem:Gfd}
  Let $\mapdef{\f}{R}{S}$ be a local homomorphism and $N$ a
  homologically finite $S$--complex.  If $\,\Gfd{N}$ is finite, then
  $\Gfd[\Rhat]{\tpP[S]{\Shat}{N}}$ is finite as well, and there is an
  inequality: $\Gfd{N} \le \Gfd[\Rhat]{\tpP[S]{\Shat}{N}}$.
  \begin{spg}
    In \corref{Gfd} we strengthen the inequality to an equality.
  \end{spg}
\end{lem}

\begin{prf}
  By \prpcite[8.13]{SInSSW04} the G--dimension of $N$ along $\f$,
  introduced in that paper and denoted $\Gdim[\f]{N}$, is finite.  By
  \rescite[3.4.1]{SInSSW04} also $\Gdim[\hat{\f}]{\tpP[S]{\Shat}{N}}$
  is finite, where $\mapdef{\hat{\f}}{\Rhat}{\Shat}$ is the completion
  of $\f$, and hence $\Gfd[\Rhat]{\tpP[S]{\Shat}{N}}$ is finite, by
  \thmcite[8.2]{SInSSW04}. Moreover, we have
  \begin{equation*}
    \Gfd{N} = \Rfd{N} \le \Rfd[\Rhat]{\tpP[S]{\Shat}{N}} =
    \Gfd[\Rhat]{\tpP[S]{\Shat}{N}},
  \end{equation*}
  where the equalities are by \thmcite[8.8]{SInSSW04} and the
  inequality is \lempartref{Rfd}{b}.
\end{prf}

\vspace*{0pt}
\begin{prf}[Proof of \thmref{main}]
  By hypothesis $N$ is an $R$--complex and there exists a local
  homomorphism $\mapdef{\f}{R}{S}$ such that $N$ is a homologically
  finite $S$--complex. It suffices to prove
  \begin{equation*}
    \tag{\one}
    \Gfd{N}  = \supP{\Dtp{\E{k}}{N}},
  \end{equation*}
  since the second equality of the claim is \thmcite[8.7]{SInSSW04}.
  
  \begin{spg}
    \step{1} First we reduce the problem to the case where $R$ and $S$
    are complete (in the topologies induced by the respective maximal
    ideals). The right hand side in (\one) is unchanged on tensoring
    with $\Shat$: Indeed there are isomorphisms of complexes
    \begin{align*}
      \tp[S]{\DtpP{\E{k}}{N}}{\Shat} &\eq
      \Dtp{\E{k}}{\tpP[S]{N}{\Shat}} \\
      &\eq \Dtp{\E{k}}{\tpP[S]{\Shat}{N}}\\
      &\eq \Dtp{\E{k}}{\DtpP[\Rhat]{\Rhat}{\tpP[S]{\Shat}{N}}}\\
      &\eq \Dtp[\Rhat]{\DtpP{\E{k}}{\Rhat}}{\tpP[S]{\Shat}{N}}\\
      &\eq \Dtp[\Rhat]{\E[\Rhat]{k}}{\tpP[S]{\Shat}{N}}
    \end{align*}
    where the first and penultimate ones hold by associativity of tensor
    products. The second isomorphism holds as $\tp[S]{N}{\Shat}$ and
    $\tp[S]{\Shat}{N}$ are isomorphic as $S$--complexes and hence also
    as $R$--complexes. The third isomorphism holds because the composite
    map $R \to S \to \Shat$ factors through $\Rhat$. Being $\m$-torsion,
    $\E{k}$ is naturally isomorphic to $\tp{\E{k}}{\Rhat}$, and as
    $\Rhat$--modules $\E{k} \is \E[\Rhat]{k}$; this accounts for the
    last isomorphism. The faithful flatness of $\Shat$ over $S$ and the
    isomorphisms above yield:
    \begin{align*}
      \supP{\Dtp{\E{k}}{N}} &= \supP{\tp[S]{\DtpP{\E{k}}{N}}{\Shat}}\\
      &= \supP{\Dtp[\Rhat]{\E[\Rhat]{k}}{\tpP[S]{\Shat}{N}}}
    \end{align*}
    The preceding equality and \rmkref{holm} yield the first two
    (in)equalities below, while \lemref{Gfd} gives the third one:
    \begin{equation*}
      \tag{\two}
      \begin{split}
        \supP{\Dtp[\Rhat]{\E[\Rhat]{k}}{\tpP[S]{\Shat}{N}}} &=
        \supP{\Dtp{\E{k}}{N}} \\ &\le \Gfd{N} \\ &\le
        \Gfd[\Rhat]{\tpP[S]{\Shat}{N}}.
      \end{split}
    \end{equation*}
    Moreover, $\Gfd[\Rhat]{\tpP[S]{\Shat}{N}}$ is finite, again by
    \lemref[]{Gfd}, and the complex $\tp[S]{\Shat}{N}$ is homologically
    finite over $\Shat$.  Thus, if (\one) holds when $R$ and $S$ are
    complete, then equalities must hold all way through in (\two).
    
    We assume henceforth that $R$ and $S$ are complete.
  \end{spg}

  \begin{spg}
    \step{2} Next we reduce to the case where $\f$ is flat and the
    closed fiber $S/\m S$ is regular. Since $R$ and $S$ are complete,
    the homomorphism $\f$ admits a regular factorization: a
    commutative diagram of local homomorphisms
    \begin{displaymath}
      \xymatrix{ {} & R' \ar@{->>}[dr]^-{\f'} & {} \\
        R \ar[ur]^-{\dot{\f}} \ar[rr]_-{\varphi} & {} & S }
    \end{displaymath}
    where $\f'$ is surjective and $\dot{\f}$ is flat with $R'/\m R'$
    regular, cf.~\thmcite[(1.1)]{AFH-94}. Since $N$ is homologically
    finite over $S$ it is also finite over $R'$, and so it suffices to
    prove the result for $\dot{\f}$; this achieves the desired
    reduction.
  \end{spg}

  \begin{spg}
    \step{3} Since $R$ is complete, it has a dualizing complex $D$;
    since $\f$ is flat with regular closed fiber, the complex
    $\tp{S}{D}$ is dualizing for $S$, cf.~\cite{rad}. Now, from
    \prpcite[(5.3)]{LWC01a} it follows that an $S$--complex $X$ is in
    the Auslander category $\A[S]$ if and only if it is in $\A$. By
    \thmcite[(4.1)]{CFH-} complexes in the Auslander category are
    exactly those of finite Gorenstein flat dimension, that is,
    \begin{align*}
      \tag{\five}
      \begin{split}
        \Gfd{X}<\infty &\Lra X\in\A \\
        &\Lra X\in\A[S] \\
        &\Lra \Gfd[S]{X}<\infty.
      \end{split}
    \end{align*}
    Therefore, when $\Gfd{X}$ is finite, so is $\Gfd[S]{X}$, and hence
    \begin{align*}
      \tag{\three} \Gfd{X} = \Rfd{X} \le \Rfd[S]{X} = \Gfd[S]{X},
    \end{align*}
    where the inequality is \lempartref{Rfd}{a} and the equalities are
    by \thmcite[8.8]{SInSSW04}.
    
    We may assume that $\H{N}\ne 0$ and set $i=\inf{N}$. By (\five)
    the complex $N$ has finite Gorenstein flat dimension over $S$;
    since it is homologically finite, it thus has finite G--dimension
    over $S$, cf.~\thmcite[(5.1.11)]{LWC}. By \thmref{fgapprox} there
    is an exact triangle in $\Df[S]$:
    \begin{align*}
      N \to P \to H \to \Susp{}{N},
    \end{align*}
    where $\pd[S]{P} = \Gdim[S]{N}$ and $\Gdim[S]{H} \le i$; in
    particular $\Gfd[S]{H}\le i$, again by \thmcite[(5.1.11)]{LWC}. By
    (\three) it follows that $\Gfd{H} \le \Gfd[S]{H} \le i$.  For any
    injective $R$--module $J$ one therefore has $\supP{\Dtp{J}{H}}\le
    i$ by \rmkref[]{holm}, and hence the exact triangle above yields
    the following isomorphisms and exact sequence
    \begin{equation*}
      \tag{\seven}
      \begin{split}
        \H[\l]{\Dtp{J}{N}} &\is \H[\l]{\Dtp{J}{P}} \quad
        \text{for}\quad \l \ge i+1, \\
        0 \to \H[i]{\Dtp{J}{N}} &\to \H[i]{\Dtp{J}{P}}.
      \end{split}
    \end{equation*}
    Since $\inf{P} \ge i$ we deduce that $\sup{\DtpP{J}{N}} \le
    \sup{\DtpP{J}{P}}$. Combined with \rmkref[]{holm} this implies the
    second inequality below
    \begin{align*}
      \supP{\Dtp{\E{k}}{N}} &\le \Gfd{N}\\
      &\le \Gfd{P}\\
      &\le \fd{P}\\
      &= \sup{\DtpP{\E{k}}{P}};
    \end{align*}
    the first inequality is also by \rmkref[]{holm}, the third
    inequality is trivial, while the equality is by \eqref{fd}, since
    $\f$ flat and $\pd[S]{P}$ finite implies $\fd{P}$ finite.
    
    Finally, $\H[\l]{\Dtp{\E{k}}{N}} \ne 0$ for some $\l \ge
    i=\inf{N}$, cf.~\lemref{supp}, and so (\seven) shows that
    $\supP{\Dtp{\E{k}}{N}} = \sup{\DtpP{\E{k}}{P}}$.  Thus, from the
    preceding display, we conclude that $\supP{\Dtp{\E{k}}{N}} =
    \Gfd{N}$.\qedhere
  \end{spg}
\end{prf}

\begin{cor}
  \label{cor:Gfd}
  Let $\mapdef{\f}{R}{S}$ be a local homomorphism and $N$ a
  homologically finite $S$--complex. If $\,\Gfd{N}$ is finite, then
  \begin{equation*}
    \Gfd{N} = \Gfd[\Rhat]{\tpP[S]{\Shat}{N}}.
  \end{equation*} 
\end{cor}

\begin{prf}
  From \lemref{Gfd} one obtains that $\Gfd[\Rhat]{\tpP[S]{\Shat}{N}}$
  is finite. Since the $\Rhat$--complex $\tp[S]{\Shat}{N}$ is finite
  over the completion $\mapdef{\hat{\f}}{\Rhat}{\Shat}$, \thmref{main}
  gives the first and the last equalities below:
  \begin{align*}
    \Gfd{N} &= \dptR - \dpt{N} \\
    &= \dpt[\no]{\Rhat} - \dpt[\Rhat]{\tpP[S]{\Shat}{N}} \\
    &= \Gfd[\Rhat]{\tpP[S]{\Shat}{N}};
  \end{align*}
  the second equality is a standard property of depth.
\end{prf}

\begin{ipg}
  We conclude with a global version of \thmref{main}:
\end{ipg}

\begin{thm}
  Let $\mapdef{\f}{R}{S}$ be a homomorphism of rings and let $N$ be a
  homologically finite $S$--complex. If $\,\Gfd{N}$ is finite, then
  \begin{equation*}
    \Gfd{X} = \supremum{\supP{\Dtp[R_\p]{\E{k(\p)}}{N_\q}}}{\q\in
      \Max[S]{N}\;\text{\rm and}\;\,\p=\q\cap R}.
  \end{equation*}
  \begin{spg}
    Note that $\msupp[S]{N} = \Max[S]{N}$ as $N$ is homologically
    finite.
  \end{spg}
\end{thm}

\begin{prf}
  For each $\q\in\Max[S]{N}$ the $R_\p$--complex $N_\q$ is finite over
  the local homomorphism $\mapdef{\f_\q}{R_\p}{S_\q}$, and
  $\Gfd[R_\p]{N_\q}$ is finite by \thmref{loc}, so \thmref[]{main}
  yields
  \[
  \Gfd[R_\p]{N_\q} =\supP{\Dtp[R_\p]{\E{k(\p)}}{N_\q}}.
  \]
  Combining this equality with that in \thmref{loc} gives the desired
  result.
\end{prf}

\section*{Acknowledgments}

We thank Lucho Avramov and Sean Sather-Wagstaff for their comments and
suggestions on this work.


\bibliographystyle{amsplain}

\providecommand{\bysame}{\leavevmode\hbox to3em{\hrulefill}\thinspace}
\providecommand{\MR}{\relax\ifhmode\unskip\space\fi MR }
\providecommand{\MRhref}[2]{%
  \href{http://www.ams.org/mathscinet-getitem?mr=#1}{#2} }
\providecommand{\href}[2]{#2}

\end{document}